\documentclass[11pt,twoside]{article}
\usepackage{latexsym}
\usepackage{amssymb,amsbsy,amsmath,amsfonts,amssymb,amscd}
\usepackage{epsfig, graphicx,subfigure}
\usepackage[active]{srcltx}
\usepackage{chemarrow}
\usepackage{psfrag}

\usepackage{colortbl}
\newcommand{\commentout}[1]{}

\newcommand{\1}{{\mathchoice {\rm 1\mskip-4mu l} {\rm 1\mskip-4mu l}
{\rm 1\mskip-4.5mu l} {\rm 1\mskip-5mu l}}}

\newcommand {\Chi} {{\bf \raise 2pt \hbox{$\chi$}} }

\newcommand {\T} { {\mathcal T} }
\newcommand {\Q} { {\mathcal Q} }
\newcommand {\F} { {\mathcal F} }
\newcommand {\C} { {\mathcal C} }
\newcommand {\G} { {\mathcal G} }
\newcommand {\f}   {\frac}
\newcommand {\p}   {\partial}

\newcommand{\dis}{\displaystyle}

\newcommand{\beq}{\begin{equation}}
\newcommand{\beqa} {\begin{array}{rl}}
\newcommand{\eeq}{\end{equation}}
\newcommand{\eeqa}{\end{array}}

\makeatletter
\newcommand{\xleftrightarrows}[2][]{\mathrel{
 \raise.40ex\hbox{$
       \ext@arrow 3095\leftarrowfill@{\phantom{#1}}{#2}$}
 \setbox0=\hbox{$\ext@arrow 0359\rightarrowfill@{#1}{\phantom{#2}}$}
 \kern-\wd0 \lower.40ex\box0}}

\newcommand{\xrightleftarrows}[2][]{\mathrel{
 \raise.40ex\hbox{$\ext@arrow 3095\rightarrowfill@{\phantom{#1}}{#2}$}
 \setbox0=\hbox{$\ext@arrow 0359\leftarrowfill@{#1}{\phantom{#2}}$}
 \kern-\wd0 \lower.40ex\box0}}
\newcommand{\xleftrightarrow}[2][]{
     \ext@arrow 0055{\leftrightarrowfill@}{#1}{#2}
}

\newcommand{\qed}{{ \hfill
                       {\unskip\kern 6pt\penalty 500
                       \raise -2pt\hbox{\vrule\vbox to 6pt{\hrule width 6pt
                       \vfill\hrule}\vrule} \par}   }}

\title{\LARGE High-order WENO scheme for Polymerization-type equations}

\author{Pierre Gabriel\thanks{Universit\'e Pierre et Marie Curie-Paris 6, UMR 7598 LJLL, BC187, 4, place de Jussieu, F-75252 Paris cedex 5, France; email: gabriel@ann.jussieu.fr} \and L\'eon Matar Tine\thanks{Equipe-Projet SIMPAF, Centre de Recherche INRIA Lille Nord Europe, Parc Scientifique de la Haute Borne, 40, avenue Halley B.P. 70478, F-59658 Villeneuve d’Ascq cedex, France; email: leonmatar@yahoo.fr} \thanks{Laboratoire Paul Painlev\'é, UMR 8524, CNRS–Universit\'é des Sciences et Technologies de Lille Cit\'é Scientifique, F-59655 Villeneuve d’Ascq Cedex, France; email: Leon-matar.Tine@math.univ-lille1.fr} \thanks{Laboratoire LANI - UFR SAT, Universit\'e Gaston Berger de Saint-Louis, B. P. 234 Saint-Louis, S\'en\'egal.}}

\date{\today}

\begin{document}

\maketitle

\pagestyle{plain}
\pagenumbering{arabic}

\begin{abstract}
Polymerization of proteins is a biochemical process involved in different diseases. Mathematically, it is generally modeled by aggregation-fragmentation-type equations. In this paper we consider a general polymerization model and propose a high-order numerical scheme to investigate the behavior of the solution. An important property of the equation is the mass conservation. The fifth-order WENO scheme is built to preserve the total mass of proteins along time.
\end{abstract}

\


\

\noindent{\bf AMS Class. No.} 35B40, 35F50, 35L65, 35M30, 35Q92, 35R09, 65-06, 65M06, 65R20

\

\section*{Introduction}

The central mechanism of amyloid diseases is the polymerization of proteins : PrP in Prion diseases, APP in Alzheimer, Htt in Huntington. The abnormal form of these proteins is pathogenic and has the ability to polymerize into fibrils. In order to well understand this process, investigation of the size repartition of polymers is a crucial point. To this end, we discuss in this paper the mathematical modeling of these polymerization processes and we propose numerical methods to investigate the mathematical features of the models.

Mathematical models are already widely used to study the polymerization mechanism of Prion diseases \cite{CL1,DGL,Engler,Greer,Greer2,LW,Masel3,Masel,Lenuzza,Pruss}, Alzheimer \cite{Craft,Lomakin,Pallitto} or Huntington \cite{Cajavec}. Such models are also used for other biological polymerization processes \cite{Banasiak, biben} and even for cell division \cite{Bekkal1, Doumic, BP} or in neurosciences \cite{PPS}.

Another field where we find aggregation-fragmentation equations is the physics of aggregates (aerosol and raindrop formation, smoke, sprays...). Among these models (see \cite{LM2} for a review), one can mention the Smoluchowsky coagulation equation \cite{EscoMischler4,FL2,FL1,LM1,MischlerRicard} with fragmentation \cite{EscoMischler1,EscoMischler2,EscoMischler3,LM4,LM3,LM5}
and the Lifshitz-Slyosov system \cite{Carrillo,CG1,Collet,FL3,NP2,NP1,NV2,NV1}. In \cite{CG2,Herrmann} a Smoluchowsky coagulation term is added to the Lifshitz-Slyosov equation.

In this paper we are interested in a model including polymerization, coagulation and fragmentation phenomena. We consider a medium where there are monomers (normal proteins for instance) characterized by the concentration $V(t)$ at time $t$ and polymers (aggregates of abnormal proteins) of size $x$ with the concentration $u(t,x).$ The dynamics of the density function $u(t,x)$ is driven by the system
\begin{equation}
\label{eq:model}
\left\{
\begin{array}{l}
\begin{array}{lll}
\dfrac{d}{dt}V(t) &=& \displaystyle -\int_{0}^{\infty} \T\bigl(V(t),x\bigr)u(t,x) \; dx,
\vspace{.4cm}\\
\dfrac{\partial}{\partial t} u(t,x) &=& -\displaystyle \frac{\partial}{\partial x} \Bigl(\T\bigl(V(t),x\bigr)u(t,x)\Bigr)+\Q(u)(t,x),\end{array}
\vspace{.7cm}\\
u(t,0)=0,\quad u(0,x)=u_0(x)\geq0\quad\text{and}\quad V(0)=V_0\geq0.
\end{array} \right.
\end{equation}
The monomers are attached by polymers of size $x$ with the polymerization rate $k_\text{on}(x).$ Depolymerization occurs when monomers detach from polymers with a rate $k_\text{off}(x).$ Hence the transport term writes
\beq\label{eq:polymerization}\T(V,x)=Vk_{\text{on}}(x)-k_{\text{off}}(x).\eeq
The two functions $k_\text{on}$ and $k_\text{off}$ are piecewise derivable but can be discontinuous. They are positive except that $k_\text{on}$ can vanish at zero. In this case, or more generally when $\T(V(t),0)\leq0,$ the boundary condition on $u(t,0)$ is not necessary since the characteristic curves outgo from the domain. The choice of the boundary condition $u(t,0)=0$ is justified later.\\
The coalescence of two polymers and the fragmentation of a polymer into two smaller ones are taken into account by the operator $\Q.$ More precisely, denoting by $A_x$ an aggregate of size $x$ we have
\begin{eqnarray*}
A_x\ +\ A_y&\xrightarrow{k_\text{c}(x,y)}&A_{x+y}\hspace{2.4cm}\text{coagulation}\\
A_{x+y}&\xrightarrow{k_\text{f}(x,y)}&A_x\ +\ A_y\hspace{1.6cm}\text{fragmentation}.
\end{eqnarray*}
Thus the coagulation-fragmentation operator is $\Q=\Q_\text{c}-\Q_\text{f}$ with
\beq\label{eq:coagulation}\Q_\text{c}(u)(x)= \frac12 \int_0^x k_\text{c}(y,x-y)\,u(y)u(x-y)\,dy - u(x)\int_0^\infty k_\text{c}(x,y)\,u(y)\,dy,\eeq
\beq\label{eq:fragmentation}\Q_\text{f}(u)(x)= \f12u(x)\int_0^x k_\text{f}(y,x-y)\,dy - \int_{0}^{\infty} k_\text{f}(x,y) \, u(x+y) \, dy.\eeq
The coalescence of two polymers of size $x$ and $y$ occurs with the symmetric rate $k_\text{c}(x,y)=k_\text{c}(y,x).$ This rate is a nonnegative function as the fragmentation symmetric rate $k_\text{f}(x,y)=k_\text{f}(y,x)$ with which a polymer of size $x+y$ produces two fragments of size $x$ and $y.$

\

There is a difference between $V(t)$ and $u(t,x=0).$ In biochemical polymerization processes, small polymers are very unstable and thus do not exist. When they appear by detachment from a longer polymer, they are immediately degraded into monomers. Thus, the quantity of small polymers vanishes while the quantity of monomers is very high. To reflect this in the mathematical model, a quantity $V(t)$ of monomers is introduced, which is different from the quantity of small polymers $u(t,x=0).$ The evolution of the first one is given by an ODE while the second one is forced to be equal to zero through the boundary condition $u(t,0)=0.$ A consequence of this distinction is that starting from $u_0(x)=0$ and $V_0>0$ there is no evolution~: the concentration of monomers is constant in time, $V(t)=V_0,$ and the concentration of polymers remains null, $u(t,x)=0.$ This is a very  intuitive and natural behaviour which is important to preserve for biological applications.\\
In the modeling, the distinction between $V$ and $u(x=0)$ induces a separation of the polymerization-depolymerization process from the coagulation-fragmentation. Indeed the aggregation of a monomer to a polymer can be seen as a coagulation but the resulting polymer has same size $x$ than the initial one, since a monomer is very small compared to the typical size of a polymer. So a transport term is more accurate to model this phenomenon than an integral term (see \cite{DGL}).\\
There is also the fact that when a small polymer is degraded into monomers, it increases the quantity of monomers. In a discrete model, this term appears in the equation on $V$ (see $n_0$ in \cite{Masel}). In the continuous model~\eqref{eq:model} this term can be neglected since the quantity of monomers produced by degradation of small polymers is very small compared to the total quantity of monomers.

\

\section{Mass Conservation}

The mechanism of polymerization is nothing but a rearrangement of the proteins, there is no creation and no disparition. So the total quantity of proteins has to be constant in time and this is the case in model~\eqref{eq:model}. We define the total mass of the system as
\beq\label{eq:totalmass}P(t)=V(t)+\int_0^\infty xu(t,x)\,dx,\eeq
since a polymer of size $x$ ``contains $x$ monomers''.
Integrating the equation on $u(t,x)$ multiplied by $x$ and adding the equation on $V$ we obtain
\beq\label{eq:masspreserving}\forall t>0,\qquad \f{dP}{dt}(t)=0,\eeq
so the total mass is conserved along time. This is a very important property that we want to keep in the numerical scheme and for this we rewrite equation~\eqref{eq:model} under a conservative form.

\subsection{Conservative formulation}

The classical discretization methods for transport equations are mass preserving. So the idea is to write the coagulation-fragmentation operator $\Q,$ which preserves the mass, under a conservative form in order to use a transport scheme. For this we follow the paper \cite{Filbet} where such a transformation is made :
\begin{equation*}
\left\{
\begin{array}{l}
\dis x\Q_{\text{c}}(u)(x)=-\f{\p\C(u)}{\p x}(x),
\vspace{.2cm}\\
\dis x\Q_{\text{f}}(u)(x)=-\f{\p\F(u)}{\p x}(x),
\end{array} \right.
\end{equation*}
where the operator $\C(u)$ is given by
\beq\label{eq:coagcons}\C(u)(x):=\int_0^x\int_{x-y}^\infty yk_\text{c}(y,z)u(y)u(z)\,dzdy,\eeq
and $\F(u)$ is
\beq\label{eq:fragcons}\F(u)(x):=\int_0^x\int_{x-y}^\infty yk_\text{f}(y,z)u(y+z)\,dzdy.\eeq
Under this form, the mass conservation is clearer and the use of conservative schemes possible.

A useful consequence of the property~\eqref{eq:masspreserving} for the numerical scheme is that the ODE on $V$ can be replaced by a mass conservation equation (see \cite{Herrmann})
\beq\label{eq:masspreserving2}\forall t>0,\qquad V(t)=V_0+\int_0^\infty x(u_0(x)-u(t,x))\,dx.\eeq
Numerically, this equation is much easier to compute than the ODE to be solved. Moreover \eqref{eq:masspreserving2} provides an explicit expression for $V$ as a function of $u.$ So we set
\beq\label{eq:polymerization2}\G(u)(x):=\Bigl(V_0+\int_0^\infty y\bigl[u_0(y)-u(y)\bigr]dy\Bigr)k_\text{on}(x)-k_\text{off}(x)\eeq
and we obtain a new equation equivalent to \eqref{eq:model}
\begin{equation}\label{eq:cons}
\left\{\begin{array}{l}
\dis x\f{\p}{\p t} u(t,x) + \f{\p\bigl[\G(u)xu\bigr]}{\p x} (t,x) + \f{\p\C(u)}{\p x}(t,x) - \f{\p\F(u)}{\p x}(t,x) = \G(u)u(t,x),
\vspace{.3cm}\\
\dis u(t,0)=0,\qquad u(0,x)=u_0(x).
\end{array}\right.
\end{equation}
In this equation~\eqref{eq:cons}, we have written the transport term as
\beq\label{eq:transport}x\f{\p\bigl[\G(u)u\bigr]}{\p x} (t,x)=\f{\p\bigl[\G(u)xu\bigr]}{\p x} (t,x)-\G(u)u(t,x).\eeq
This formulation enhances the relation
\beq\label{eq:masspreserving3}\f{d}{dt}\int_0^\infty xu(t,x)\,dx=\int_0^\infty \G(u)u(t,x)\,dx=-\f{d}{dt}V(t)\eeq
and allows to preserve this property numerically when using conservative transport schemes.

\subsection{Domain truncation}

Numerically, equation~\eqref{eq:cons} is solved on a truncated domain $[0,R]$ so the integration bounds have to be changed in order to keep the mass preservation. For the coagulation term, we introduce as in \cite{Filbet}
\begin{eqnarray*}
\C^R(u)(x)&:=&\dis\int_0^x\int_{x-y}^{R-y} yk_\text{c}(y,z)u(y)u(z)\,dzdy\\
&=&\dis\int_0^x\int_{x}^{R} yk_\text{c}(y,z-y)u(y)u(z-y)\,dzdy,
\end{eqnarray*}
and for the fragmentation
\begin{eqnarray*}
\F^R(u)(x)&:=&\dis\int_0^x\int_{x-y}^{R-y} yk_\text{f}(y,z)u(y+z)\,dzdy\\
&=&\dis\int_0^x\int_{x}^{R} yk_\text{f}(y,z-y)u(z)\,dzdy.
\end{eqnarray*}
With this truncation, we have $\C^R(u)(0)=\C^R(u)(R)=\F^R(u)(0)=\F^R(u)(R)=0.$ So the total mass does neither increase nor decrease with respect to time if we consider the coagulation and fragmentation processes. If we look at the effects of this truncation on the original coagulation and fragmentation operators we have
$$Q^R_\text{c}(u)(x):=-\f1x\p_x\C^R(u)(x)=\f12\int_0^x k_\text{c}(y,x-y)u(u)u(x-y)\,dy-u(x)\int_0^{R-x}k_\text{c}(x,y)u(y)\,dy$$
and
$$Q^R_\text{f}(u)(x):=-\f1x\p_x\F^R(u)(x)=\f12u(x)\int_0^x k_\text{f}(y,x-y)\,dy-\int_x^{R}k_\text{f}(x,y-x)u(y)\,dy.$$
In the coagulation term, the truncation corresponds to the assumption that a polymer of size $x$ cannot coagulate with a polymer of size greater than $R-x.$ Concerning the fragmentation term, it is nothing but the assumption that polymers of size greater than $R$ cannot split. Biologically they are the natural assumptions to avoid the loss of mass.

Concerning the transport term, the only way to avoid the loss of mass is to set
\beq\label{eq:GxR}\G^R(u)(R,t)=0.\eeq
The meaning we give to this relation in the numerical scheme is exposed in Section~\ref{sse:WENO}. It is useless to do such a truncation for $x=0$ since $xu(t,x)$ vanishes when $x=0.$

Finally we obtain a conservative truncated equation for $x\in(0,R)$
\begin{equation}\label{eq:truncated}
\left\{\begin{array}{l}
\dis x\f{\p}{\p t} u_R(t,x) + \f{\p\bigl[\G^R(u_R)xu_R\bigr]}{\p x} (t,x) + \f{\p\C^R(u_R)}{\p x}(t,x) - \f{\p\F^R(u_R)}{\p x}(t,x) = \G(u_R)u_R(t,x),
\vspace{.3cm}\\
\dis u_R(t,0)=0,\qquad u_R(0,x)=u_0(x).
\end{array}\right.
\end{equation}

\

When there is no transport term, convergence of the solution of Equation~\eqref{eq:truncated} to the solution of Equation~\eqref{eq:model} when $R\to\infty$ is proved in \cite{DubovStew,LM2,LM4,LM5,Stewart} under growth conditions on $k_\text{c}$ and $k_\text{f}.$

\

\section{A High Order WENO Scheme}\label{sec:WENO}

In order to obtain a mass preserving scheme, we consider equation~\eqref{eq:cons} as a transport equation and for high order accuracy we choose a fifth-order WENO (Weighted Essentially Non Oscillatory) reconstruction for the fluxes. This high order scheme is comonly used \cite{Devys,SSZ} since it is not more complicated to implement than a third order WENO one for instance.

\subsection{Numerical fluxes}

Before using the WENO reconstruction we have to know if the fluxes are positive or negative in order to appropriately upwind the scheme. Concerning the coagulation and the fragmentation terms, we consider a positive upwinding as suggested in \cite{Filbet}. For the transport term $\p_x\bigl[\G(u)xu\bigr]$ we have to make a flux splitting because $\G$ has no sign. A natural splitting here is to put the terms of $\G$ that are preceded by a plus sign in the positive part and the terms preceded by a minus sign in the negative part, namely $\G=\G_1^++\G_1^-$ where
\beq\label{eq:flux1}\left\{\begin{array}{l}
\G_1^+(u)(x)=\Bigl(V_0+\int_0^\infty yu_0(y)dy\Bigr)k_\text{on}(x),\\
\G_1^-(u)(x)=-\Bigl(\int_0^\infty yu(y)dy\Bigr)k_\text{on}(x)-k_\text{off}(x).
\end{array}\right.\eeq
An other decomposition is the polymerization-depolymerization one
\beq\label{eq:flux0}\left\{\begin{array}{l}
\G_0^+(u)(x)=\Bigl(V_0+\int_0^\infty y\bigl[u_0(y)-u(y)\bigr]dy\Bigr)k_\text{on}(x),
\vspace{.15cm}\\
\G_0^-(u)(x)=-k_\text{off}(x).
\end{array}\right.\eeq
The term $\G_0^+$ is necessarily positive because $V_0+\int_0^\infty y\bigl[u_0(y)-u(t,y)\bigr]dy=V(t)\geq0.$ With these two flux splittings, we built others by convex combination. For any $\lambda\in[0,1]$ we set $\G_\lambda=\lambda\G_1+(1-\lambda)\G_0$ which gives
\beq\label{eq:fluxlb}\left\{\begin{array}{l}
\G_\lambda^+(u)(x)=\Bigl(V_0+\int_0^\infty y\bigl[u_0(y)-u(y)\bigr]dy+\lambda\int_0^\infty yu(y)dy\Bigr)k_\text{on}(x),\\
\G_\lambda^-(u)(x)=-\lambda\Bigl(\int_0^\infty yu(y)dy\Bigr)k_\text{on}(x)-k_\text{off}(x).
\end{array}\right.\eeq
We also consider the Lax-Friedrichs scheme which corresponds to
\beq\label{eq:LF}\left\{\begin{array}{l}
\G_\text{\tiny LF}^+(u)=\f12\bigl(\G(u)+m\bigr),
\vspace{.32cm}\\
\G_\text{\tiny LF}^-(u)=\f12\bigl(\G(u)-m\bigl),
\end{array}\right.\eeq
with $m=\max_{x\geq0}|\G(u)|.$ This term has to be computed at each time step because $\G(u)$ depends on time.

\

Finally, the WENO reconstruction is done with the fluxes
\beq\label{eq:fluxes}\left\{\begin{array}{l}
H^+(u)=\G^+(u)xu+\C(u)-F(u),
\vspace{.33cm}\\
H^-(u)=\G^-(u)xu,
\end{array}\right.\eeq
and the choice among the different flux splittings is discussed in Section~\ref{sse:fluxes}.

\subsection{WENO reconstruction}\label{sse:WENO}

\begin{figure}[htbp]
\psfrag{H- left}[l]{\color{red}$H^-_{i-\f12}$}
\psfrag{H- right}[l]{\color{blue}$H^-_{i+\f12}$}
\psfrag{H+ left}[l]{\color{red}$H^+_{i-\f12}$}
\psfrag{H+ right}[l]{\color{blue}$H^+_{i+\f12}$}
\psfrag{i-3}[l]{i-3}
\psfrag{i-2}[l]{i-2}
\psfrag{i-1}[l]{i-1}
\psfrag{i}[l]{i}
\psfrag{i+1}[l]{i+1}
\psfrag{i+2}[l]{i+2}
\psfrag{i+3}[l]{i+3}
\centering
\resizebox{.7\textwidth}{!}{\includegraphics[width=.8\linewidth, height=5cm]{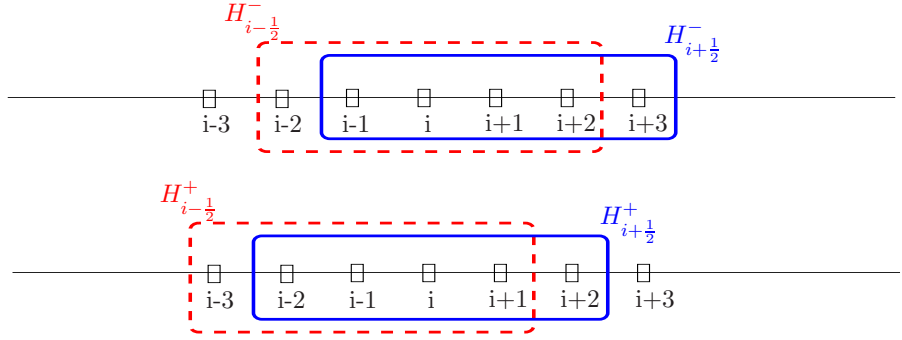}}
\caption{ stencil choice  }
\label{stencil}
\end{figure}

The point of view adopted here is the finite difference one, as recommanded in \cite{Shu}, because it is better than the finite volume in terms of operation counts. We assume the spatial domain $[0,R]$ is divided into $N$ uniform cells and we denote $x_i=i\Delta x$ for $0\leq i\leq N$ with $\Delta x=\f RN.$ We use the WENO formulation of Jiang and Peng \cite{Jiang} which consists in applying WENO to approach the spacial derivative directly on the nodes of the grid. The spatial derivative  $\p_x(H^+(v)+H^-(v))$ is approximated at the point $x_i$ by
$$\f1{\Delta x}\Bigl[H^+_{i+\f12}+H^-_{i+\f12}-H^+_{i-\f12}-H^-_{i-\f12}\Bigr]$$
where the fifth order accurate numerical flux $H^+_{i+\f12}$ is given by the WENO reconstruction. For each node $x_i$ we denote by $H^+_i$ the numerical approximation of $H^+(v(x_i)).$ The stencil choice for each flux is specified in Figure~\ref{stencil}, and the fluxes $H^{\pm}_{i\pm\f12}$ are expressed as convex combination of the $H^\pm_k$ of the stencil. Let us detail how we proceed :
\begin{itemize}
 \item [-]for $H^-_{i-\f12}$ we set $ W_{1}=H^{-}_k,\ W_{2}=H^{-}_{k+1},\ W_{3}=H^{-}_{k+2},\ W_{4}=H^{-}_{k-1},\ W_{5}=H^{-}_{k-2},$
 \item [-]for $H^-_{i+\f12}$ we set $ W_{1}=H^{-}_{k+1},\ W_{2}=H^{-}_{k+2},\ W_{3}=H^{-}_{k+3},\ W_{4}=H^{-}_k,\ W_{5}=H^{-}_{k-1},$
 \item [-]for $H^+_{i-\f12}$ we set $ W_{1}=H^{+}_{k-3},\ W_{2}=H^{+}_{k-2},\ W_{3}=H^{+}_{k-1},\ W_{4}=H^{+}_k,\ W_{5}=H^{+}_{k+1},$
 \item [-]for $H^+_{i+\f12}$ we set $ W_{1}=H^{+}_{k-2},\ W_{2}=H^{+}_{k-1},\ W_{3}=H^{+}_k,\ W_{4}=H^{+}_{k+1},\ W_{5}=H^{+}_{k+2}.$
\end{itemize}

\

For the regularity coefficients we define for each previous flux
\begin{eqnarray*}
S_{1}&=&\frac{13}{12}(W_{1}-2W_{2}+W_{3})^{2}+\frac{1}{4}(W_{1}-4W_{2}+3W_{3})^{2},\\
S_{2}&=&\frac{13}{12}(W_{2}-2W_{3}+W_{4})^{2}+\frac{1}{4}(W_{2}-W_{4})^{2},\\
S_{3}&=&\frac{13}{12}(W_{3}-2W_{4}+W_{5})^{2}+\frac{1}{4}(3W_{3}-4W_{4}+W_{5})^{2}.
\end{eqnarray*}
Then we take the weights
$$w_{r}=\frac{a_{r}}{\sum_{j=1}^{3} a_{j} }, \quad \text{with} \quad a_{r}=\frac{d_{r}}{(\epsilon+S_{r})^2}, \quad d_{1}=\frac{3}{10},\ d_{2}=\frac{6}{10},\ d_{3}=\frac{1}{10} \qquad r=1,2,3$$
where $\epsilon$ is introduced to prevent the denominator from vanishing. Finally we take the different flux parts given by
\begin{equation*}\left\{\begin{array}{l}
\dis H^-_{i\pm\f12}=w_{1}\biggl(\frac{W_{3}}{3}-\frac{7W_{2}}{6}+\frac{11W_{1}}{6}\biggr)+w_{2}\biggl(\frac{-W_{2}}{6}+\frac{5W_{1}}{6}+\frac{W_{4}}{3}\biggr)+w_{3}\biggl(\frac{W_{1}}{3}+\frac{5W_{4}}{6}-\frac{W_{5}}{6}\biggr),\\
\\
\dis H^+_{i\pm\f12}=w_{1}\biggl(\frac{W_{1}}{3}-\frac{7W_{2}}{6}+\frac{11W_{3}}{6}\biggr)+w_{2}\biggl(\frac{-W_{2}}{6}+\frac{5W_{3}}{6}+\frac{W_{4}}{3}\biggr)+w_{3}\biggl(\frac{W_{3}}{3}+\frac{5W_{4}}{6}-\frac{W_{5}}{6}\biggr).
\end{array}\right.\end{equation*}

Concerning the boundaries $x=0$ and $x=R,$ we compute the fluxes using the WENO reconstruction with ghost points $x_{-3},\ x_{-2},\ x_{-1},$ and $x_{N+1},\ x_{N+2},\ x_{N+3}.$ In the first three points we use that for all time $t\geq0,$
$$xu(t,x)_{\!^{\big|}x=0}=\C^R(u)(x=0,t)=\F^R(u)(x=0,t)=0$$
to set $H^{\pm}_{-3}=H^{\pm}_{-2}=H^{\pm}_{-1}=0.$ For the last three ones we use the truncation
$$\G^R(u)(x=R,t)=\C^R(x=R,t)=\F^R(x=R,t)=0$$
to put $H^\pm_{N+1}=H^\pm_{N+2}=H^\pm_{N+3}=0.$

\subsection{Integration method}

For the integral terms, we use a fifth order composite rule introduced in \cite{SSZ}. If $f_k$ denotes an approximation of $f(x_k),$ the method can be written as
$$\int_{i\Delta x}^{j\Delta x} f(x)\,dx\simeq \Delta x\sum^j_{k=i}\,'\, f_k$$
where
\begin{eqnarray*}
\sum^j_{k=i}\,'\, f_k&=&\f{251}{720}f_i+\f{299}{240}f_{i+1}+\f{211}{240}f_{i+2}+\f{739}{720}f_{i+3}\\
&&+\f{739}{720}f_{j-3}+\f{211}{240}f_{j-2}+\f{299}{240}f_{j-1}+\f{251}{720}f_j+\sum_{k=i+4}^{j-4}f_k
\end{eqnarray*}
if $j-i>6.$ This method is based on polynomial interpolations of the function $f.$

On the first interval, we integrate without using the boundary value $f_0=0$ because the solution can be discontinuous at $x=0.$ So we use the fifth accurate approximation
$$\sum^1_{k=0}\,'\, f_k=\f{55}{24}f_1-\f{59}{24}f_2+\f{37}{24}f_3-\f{9}{24}f_4.$$
Finally for the intervals at the boundaries we have
$$\sum^2_{0}\,'\, f_k=\f{8}{3}f_1-\f{5}{3}f_2+\f{4}{3}f_3-\f{1}{3}f_4,$$
$$\sum^3_{0}\,'\, f_k=\f{21}{8}f_1-\f{9}{8}f_2+\f{15}{8}f_3-\f{3}{8}f_4,$$
$$\sum^4_{0}\,'\, f_k=\f{21}{8}f_1-\f{7}{6}f_2+\f{29}{12}f_3+\f{1}{6}f_4-\f{1}{24}f_5,$$
$$\sum^5_{0}\,'\, f_k=\f{21}{8}f_1-\f{7}{6}f_2+\f{19}{8}f_3+\f{17}{24}f_4+\f{1}{2}f_5-\f{1}{24}f_6,$$
$$\sum^6_{0}\,'\, f_k=\f{21}{8}f_1-\f{7}{6}f_2+\f{19}{8}f_3+\f{2}{3}f_4+\f{25}{24}f_5+\f{1}{2}f_6-\f{1}{24}f_7,$$
$$\sum^7_{0}\,'\, f_k=\f{21}{8}f_1-\f{7}{6}f_2+\f{19}{8}f_3+\f{2}{3}f_4+f_5+\f{25}{24}f_6+\f{1}{2}f_7-\f{1}{24}f_8,$$
$$\sum^N_{N-1}\,'\, f_k=\f{9}{4}f_N+\f{19}{24}f_{N-1}-\f{5}{24}f_{N-2}+\f{1}{24}f_{N-3},$$
$$\sum^N_{N-2}\,'\, f_k=\f{1}{3}f_N+\f{4}{3}f_{N-1}+\f{1}{3}f_{N-2},$$
$$\sum^N_{N-3}\,'\, f_k=\f{1}{3}f_N+\f{31}{24}f_{N-1}+\f{7}{8}f_{N-2}+\f{13}{24}f_{N-3}-\f{1}{24}f_{N-4},$$
$$\sum^N_{N-4}\,'\, f_k=\f{1}{3}f_N+\f{31}{24}f_{N-1}+\f{5}{6}f_{N-2}+\f{13}{12}f_{N-3}+\f{1}{2}f_{N-4}-\f{1}{24}f_{N-5},$$
$$\sum^N_{N-5}\,'\, f_k=\f{1}{3}f_N+\f{31}{24}f_{N-1}+\f{5}{6}f_{N-2}+\f{25}{24}f_{N-3}+\f{25}{24}f_{N-4}+\f{1}{2}f_{N-5}-\f{1}{24}f_{N-6},$$
$$\sum^N_{N-6}\,'\, f_k=\f{1}{3}f_N+\f{31}{24}f_{N-1}+\f{5}{6}f_{N-2}+\f{25}{24}f_{N-3}+f_{N-4}+\f{25}{24}f_{N-5}+\f{1}{2}f_{N-6}-\f{1}{24}f_{N-7}.$$

\

We use this quadrature method to discretize the operators $\C^R$ and $\F^R$ with
\beq\label{eq:discoag}\C_i^R-\F_i^R=(\Delta x)^2\sum^i_{j=0}\,' \sum^N_{l=i+1}\,'\ x_l\bigl(k^\text{c}_{j,l-j}u_ju_{l-j}-k^\text{f}_{j,l-j}u_l\bigr)\eeq
as suggested in \cite{Filbet}. Grouping the two terms in an unique summation is lighter regarding to operation counts.

\subsection{Time discretization}

The time step is denoted by $\Delta t$ and changes along time because of the CFL stability condition that is time dependent. For the time discretization we choose a third order Runge-Kutta method. To approach the time evolution of an equation $\partial_{t} u = L(u),$ we compute at time $n\Delta t$
$$ u^{(1)}=u^{n}+\Delta t\,L(u^{n})\qquad\text{and}\qquad u^{(2)}=\frac{3}{4}u^{n}+\frac{1}{4}u^{(1)}+\frac{1}{4} \Delta t\,L(u^{(1)}),$$
where $u^n$ is an approximation of $u(n\Delta t).$ Then the approximation of $v$ at time $(n+1)\Delta t$ is given by
$$u^{n+1}= \frac{1}{3}u^{n}+\frac{2}{3}u^{(2)}+\frac{2}{3} \Delta t\,L(u^{(2)}).$$
This method is an explicit one, so to ensure the stability we compute the time step $\Delta t$ at each iteration thanks to the CFL condition
\beq\label{eq:cfl}\Delta t\leq\min{\bigl\{ (G+C+F)^{-1} \bigr\}}\eeq
where
$$G=\f1{\Delta x}\sup_{i}{\bigl(\G^+_i-\G^-_i\bigr)},\quad C=\sup_{i}{\biggl\{\sum^N_{j=1}\,'\ k^\text{c}_{i,j}u_j\biggr\}}\quad\text{and}\quad F=\sup_{i}{\biggl\{\f12\sum_{j=1}^{i-1}\,'\ k^\text{f}_{j,i-j}\biggr\}}.$$
For instance the Lax-Friedrichs decomposition leads to $G_\text{\tiny LF}=m/\Delta x.$

\

Since we combine a fifth order WENO reconstruction and a third order time discretization, we predict that our scheme is convergent of third order. To validate it numerically, we compute the solution for different discretization grids with regular parameters and initial data. Comparing these solutions at time $T=20$ in the $L^\infty$ space norm (see Table~\ref{tab:order} for the results), we obtain the numerical order $2.95$ which validates the prediction.

\begin{table}[h]
\begin{center}
\begin{tabular}{c||ccc}
\hline
$\Delta x$ & 5/40 & 5/80 & 5/160\\
\hline
error & 81.84 & 12.44 & 1.37\\
\hline
\end{tabular}

\

\caption{\label{tab:order}Error between different discrete solutions and the reference computed with $\Delta x=5/320.$}
\end{center}
\end{table}

\section{Numerical Simulations}\label{se:simul}

\subsection{Parameters}\label{sse:parameters}

Numerical values for the polymerization and fragmentation rates can be found in the biological literature (see \cite{Masel3,Rubenstein,Knowles} for instance). It is of importance to note that the models considered in these papers are discrete ones, so some computations are necessary to deduce adimensional numerical values for the continuous model~\eqref{eq:model}. Another point is that the parameters of these models do not depend on the size of polymers, so we can only obtain mean values for the size-dependent parameters.

We choose to use the values of the recent paper \cite{Knowles} to do numerical simulations. The mean length of polymers for the initial distribution is estimated to be $1380.$ With the continuous model we reduce this value to $0.2$ by considering an initial profile equal to a positive constant on $[0,0.4]$ and null for $x>0.4$ (see the first plot of Figure~\ref{fig:comparison}). Thus we define a parameter $\varepsilon=0.2/1380\approx1,4\times10^{-4}$ which allows to go from a discrete model to a continuous one (see \cite{DGL} for more details). The values we find are for instance $2.9\times10^{-2}\mu M^{-1}s^{-1}$ for the polymerization rate and $2.1\times10^{-9}s^{-1}$ for the fragmentation (where $M$ represents the concentration in $mol$ and $s$ the time in $second$). The polymerization rate appears in a derivative term, so the value of the discrete model has to be multiplied by $\varepsilon$ to obtain the continuous accurate value $4\times10^{-6}\mu M^{-1}s^{-1}.$ Conversely, the fragmentation rate which appears in an integral term has to be divided by $\varepsilon$ and we find $1.5\times10^{-5}s^{-1}.$ Concerning the depolymerization and coagulation, they are neglected in the models of \cite{Masel3,Rubenstein,Knowles}. So we consider numerical values that seem to be reasonable compared to the previous ones.

In the present study, the parameters are assumed to be size dependent as suggested in \cite{CL2,CL1} and their choice is now presented and motivated. Concerning the numerical coefficients, they are chosen in order to have mean values of the same order than the values previously obtained from \cite{Knowles}.\\
For the polymerization we assume that small polymers have a different behavior compared to the big ones. We consider a critical size $x_c=0.5$ such that polymers of size $x<x_c$ convert monomers with the rate
$$k^{(1)}_\text{on}(x)=(4x+0.2)\times10^{-6}\mu M^{-1}s^{-1},$$
and for $x>x_c$ with a constant rate
$$k^{(2)}_\text{on}(x)\equiv 4\times10^{-6}\mu M^{-1}s^{-1}.$$
For the fragmentation kernel, we use the classical assumption that the fragmentation probability depends only on the size $x+y$ of the polymer and we set
$$k_\text{f}(x,y)=\f{80(x+y)}{10+(x+y)}\times10^{-5}s^{-1}.$$
The depolymerization is assumed to be constant and of the same order as the fragmentation. We discuss the dependence on this parameter by considering different intensities
\beq\label{eq:depo}k_\text{off}(x)\equiv \eta\times10^{-6}s^{-1}\qquad\text{with}\ 2\leq \eta\leq8.\eeq
Concerning the coagulation kernel, we do not use a classical one. Even if there is no space in model~\eqref{eq:model}, we use a kernel which reflect some space effects. Indeed we consider that small polymers are very mobile and that the big ones, plaques, are very attractant. So the coagulation occurs preferentially between big and small aggregates. The kernel we choose is of the form
$$k_\text{c}(x,y)=\f{4|x-y|^\f32}{1+(x+y)}\times10^{-6}\mu M^{-1}s^{-1}.$$
This kernel satisfies the growth assumption that we can find in \cite{DubovStew} to ensure the convergence of the solution when $R\to\infty$ if we consider the only coagulation-fragmentation process.

\begin{figure}[h]

\begin{center}
\begin{minipage}{.49\textwidth}
\psfrag{'coag'}[r]{coagulation (in $\mu M^{-1}h^{-1}$)}
\psfrag{x}[l]{$x$}
\psfrag{y}[l]{$y$}
\includegraphics[width=\textwidth]{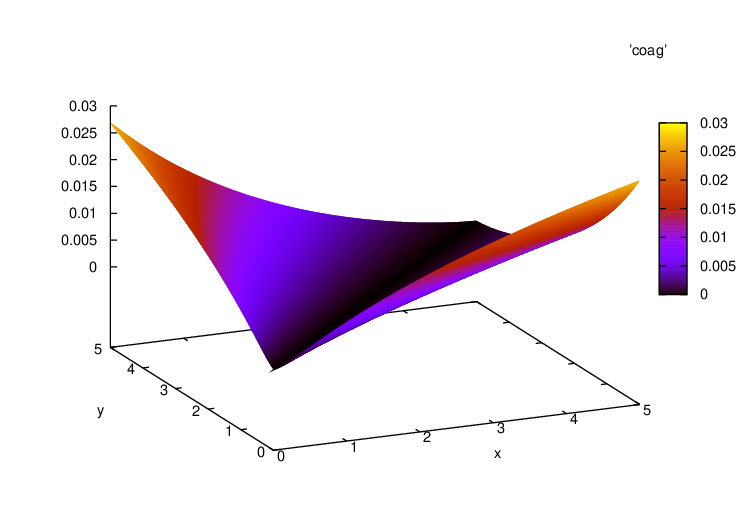}
\end{minipage}
\hfill
\begin{minipage}{.49\textwidth}
\psfrag{'frag'}[r]{fragmentation (in $h^{-1}$)}
\psfrag{x}[l]{$x$}
\psfrag{y}[l]{$y$}
\includegraphics[width=\textwidth]{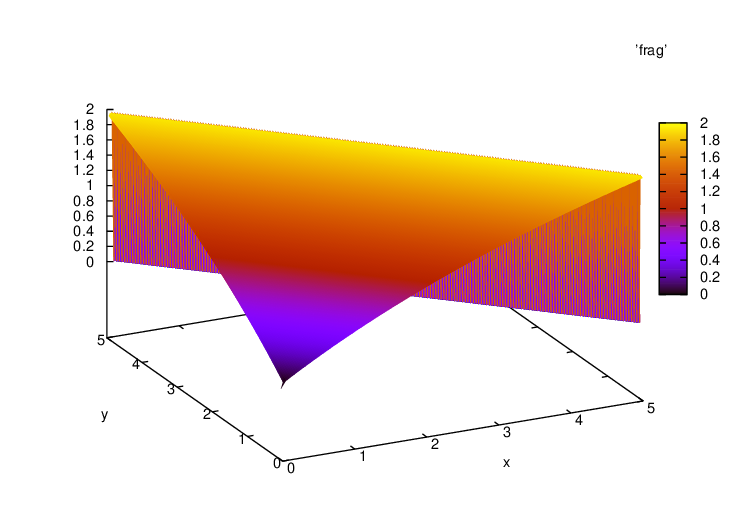}
\end{minipage}
\end{center}

\caption{Profiles of the coagulation and fragmentation kernels.}

\end{figure}

In \cite{Knowles} we also find numerical values for the initial data $V_0=98\,\mu M$ and $\dis\int_0^R xu_0(x)\,dx=0.21\mu M.$ This last value and the fact that the initial distribution of polymers is assumed to be under the form $u_0=cst\times\1_{[0,0.4]}$ lead to
\beq\label{eq:initial} u_0(x)=\left\{\begin{array}{lll}
          2.6&\text{if}&0\leq x\leq 0.4\\
          0&\text{if}&x>0.4.
         \end{array}\right.
\eeq

For the following simulations, the discretization is made on a domain $[0,5]$ with a number of nodes $N=200,$ so the mesh size is $\Delta x=0.025.$

\subsection{Choice among the different flux splittings}\label{sse:fluxes}

First we deal with the CFL condition. Thanks to the triangular inequality, we obtain that $G_\text{\tiny LF}\leq G_0.$ Moreover, there is an explicit expression for $G_\lambda$
$$G_\lambda^n=\f{1}{\Delta x}\sup_i{\biggl\{\Bigl(V_0+\Delta x\sum^N_{j=1}\,'\, x_ju_j^0+(2\lambda-1)\Delta x\sum^N_{j=1}\,'\, x_ju_j^n\Bigr)k_i^\text{on}+k_i^\text{off}\biggr\}}$$
which shows that $G_\lambda$ increases with $\lambda.$ So if $0\leq\lambda<\Lambda\leq1$ then at each time step we have \mbox{$G_\text{\tiny LF}^n\leq G_\lambda^n\leq G_\Lambda^n.$}
Notice also that, with the numerical values we have chosen, the quantity of polymers $\Delta x\sum^N_{j=1}\!'\, x_ju_j^n$ increases with $n.$ Indeed we can see in Figure~\ref{fig:monomers} that the quantity of monomers $V^n\simeq V(n\Delta t)$ defined by the mass conservation $V^n+\Delta x\sum^N_{j=1}\!'\, x_ju_j^n=V_0+\Delta x\sum^N_{j=1}\!'\, x_ju_j^0$ decreases. The consequence on the CFL condition is that $G_\lambda^n$ increases with $n$ if $\lambda>\f12,$ decreases if $\lambda<\f12$ and is time independent when $\lambda=\f12.$ Thus, regarding to the numerical computation, the fastest scheme is the Lax-Friedrichs one and then the computation time increases significantly with $\lambda.$

\

\begin{figure}[h]

\begin{minipage}{0.49\textwidth}
\psfrag{t6}[c]{$t=6$}
\psfrag{u}[c]{$u(x,6)$}
\psfrag{x}[c]{$x$}
\psfrag{L}[l]{\small $\G_\text{\tiny LF}$}
\psfrag{0}[l]{\small $\G_0$}
\psfrag{1}[l]{\small $\G_1$}
\includegraphics[width=\textwidth]{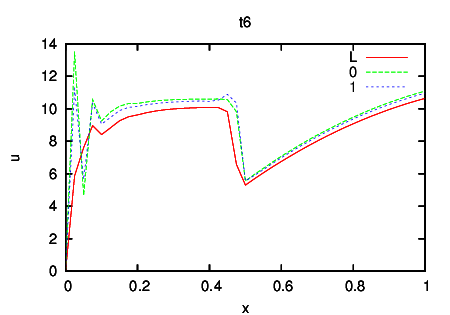}
\end{minipage}\hfill
\begin{minipage}{0.49\textwidth}
\psfrag{t12}[c]{$t=12$}
\psfrag{u}[c]{$u(x,12)$}
\psfrag{x}[c]{$x$}
\psfrag{L}[l]{\small $\G_\text{\tiny LF}$}
\psfrag{0}[l]{\small $\G_0$}
\psfrag{1}[l]{\small $\G_1$}
\includegraphics[width=\textwidth]{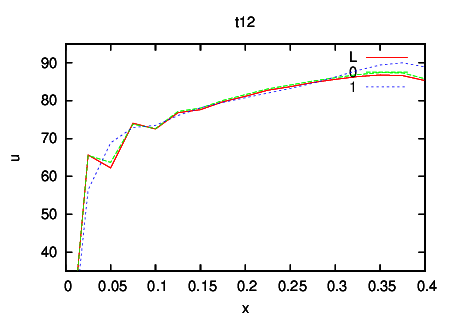}
\end{minipage}

\

\caption{Comparison between the flux splittings $\G_\text{\tiny LF},\ \G_0$ and $\G_1$ for $\eta=8\times10^{-6}s^{-1}.$}\label{fig:eta8}

\end{figure}

Let us now turn to the effects of the flux splitting on the size distribution. First we consider a depolymerization corresponding to $\eta=8$ in \eqref{eq:depo} and investigate the differences between the solutions associated to the decompositions $\G_\text{\tiny LF},\ \G_0$ and $\G_1.$ We can see in Figure~\ref{fig:eta8} that the solutions for $\G_0$ and $\G_1$ are close together for small times and then the behavior of $\G_0$ becomes closer to the Lax-Friedrichs one. The less oscillating scheme for $t=6h$ is the the Lax-Friedrichs one, but it is also the most oscillating at time $t=12h.$ Finally the solutions associated to the three flux decompositions are quite similar and they all present oscillations at some times, so we do not find with this simulation any reason to discard one of them.

\begin{figure}[h]

\begin{minipage}{0.49\textwidth}
\psfrag{t20}[c]{$\eta=6,\ t=20$}
\psfrag{u}[c]{$u(x,20)$}
\psfrag{x}[c]{$x$}
\psfrag{L}[l]{\small $\G_\text{\tiny LF}$}
\psfrag{0}[l]{\small $\G_0$}
\includegraphics[width=\textwidth]{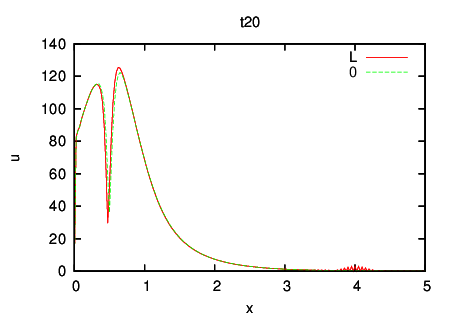}
\end{minipage}\hfill
\begin{minipage}{0.49\textwidth}
\psfrag{t18}[c]{$\eta=2,\ t=18$}
\psfrag{u}[c]{$u(x,18)$}
\psfrag{x}[c]{$x$}
\psfrag{0}[l]{\small $\G_0$}
\psfrag{2}[l]{\small $\G_{0.2}$}
\includegraphics[width=\textwidth]{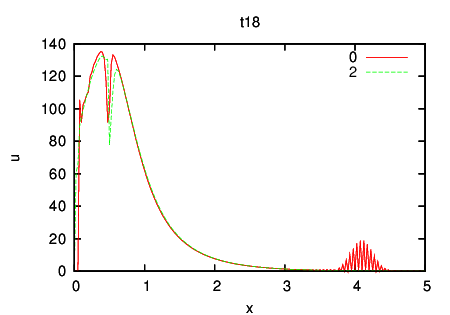}
\end{minipage}

\

\caption{Unstability of some schemes when $\eta$ decreases. Left: $\G_\text{\tiny LF}$ becomes unstable for $\eta=6.$ Right: $\G_0$ becomes unstable for $\eta=2.$}\label{fig:unstability}

\end{figure}

If we change the depolymerization rate by considering $k_\text{off}=6\times10^{-6}s^{-1},$ we remark that the Lax-Friedrichs scheme is unstable (see Figure~\ref{fig:unstability}) while $\G_0$ is stable. If we continue to decrease $\eta,$ we find with $k_\text{off}=2\times10^{-6}s^{-1}$ that the $\G_0$-scheme becomes unstable while $\G_{0.2}$ is stable. Thus the Lax-Friedrichs scheme and the $\G_\lambda$-schemes with $\lambda$ small has to be avoided to ensure stability when small depolymerization values are considered.

\

\begin{figure}[h]

\begin{minipage}{0.49\textwidth}
\psfrag{t6}[c]{$t=6$}
\psfrag{u}[c]{$u(x,6)$}
\psfrag{x}[c]{$x$}
\psfrag{2}[l]{\small $\G_{0.2}$}
\psfrag{5}[l]{\small $\G_{0.5}$}
\psfrag{1}[l]{\small $\G_1$}
\includegraphics[width=\textwidth]{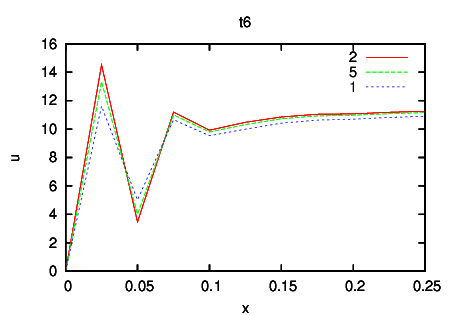}
\end{minipage}\hfill
\begin{minipage}{0.49\textwidth}
\psfrag{t20}[c]{$t=20$}
\psfrag{u}[c]{$u(x,20)$}
\psfrag{x}[c]{$x$}
\psfrag{2}[l]{\small $\G_{0.2}$}
\psfrag{5}[l]{\small $\G_{0.5}$}
\psfrag{1}[l]{\small $\G_1$}
\includegraphics[width=\textwidth]{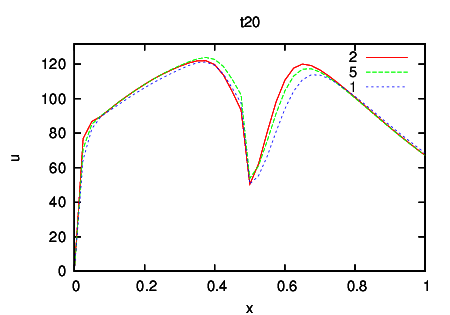}
\end{minipage}

\

\caption{Comparison of the behavior of the solution for different $\lambda$ with $\eta=5.$}\label{fig:comparison}

\end{figure}

Knowing that, we compare different stable schemes, namely $\G_\lambda$ with $0.2\leq\lambda\leq1.$ We can see in Figure~\ref{fig:comparison} that for large times ($t=20h$), the three flux decompositions provide a good behavior where there are strong variations of the solution. These locations are $x=0$ because of the boundary condition which enforces $u(t,0)$ to vanish, and $x=0.5$ because the transport term $k_\text{on}$ is discontinuous at $x=0.5.$ If we look at smaller times ($t=6h$ for instance) we can see that the larger $\lambda$ is, the less oscillating the curves are. But, as we already remarked, the quantity $G_\lambda$ is higher for $\lambda$ close to $1$ and it increases with time when $\lambda>\f12.$ Thus it is penalizing for the computation time to use high values of $\lambda.$ A good compromise could be to choose $\lambda=\f12$ since $G_\f12$ does not depend on time. The other solution is to adapt the $\lambda$ when we change the parameters.

\subsection{Interpretation of the numerical results}\label{sse:results}

\begin{figure}[h]

\begin{minipage}{0.49\textwidth}
\psfrag{t0}[c]{$t=0$}
\psfrag{u}[c]{$u(0,x)$}
\psfrag{x}[c]{$x$}
\includegraphics[width=\textwidth]{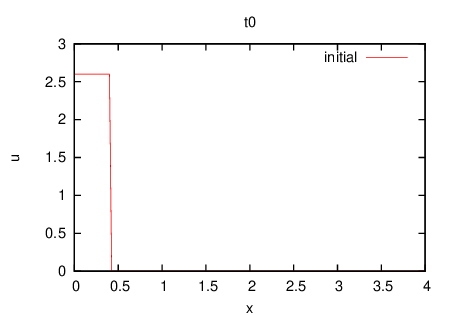}
\end{minipage}\hfill
\begin{minipage}{0.49\textwidth}
\psfrag{t1}[c]{$t=0.5$}
\psfrag{u}[c]{$u(x,0.5)$}
\psfrag{x}[c]{$x$}
\includegraphics[width=\textwidth]{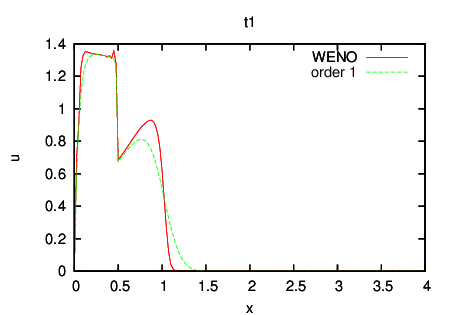}
\end{minipage}\\


\begin{minipage}{0.49\textwidth}
\psfrag{t12}[c]{$t=12$}
\psfrag{u}[c]{$u(x,12)$}
\psfrag{x}[c]{$x$}
\hspace{.1cm}\includegraphics[width=\textwidth]{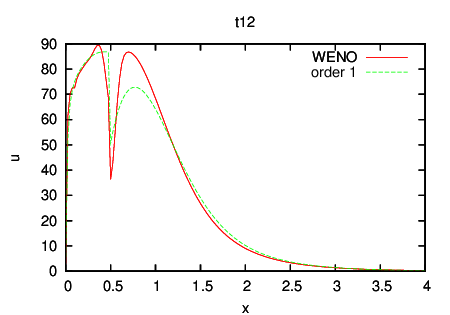}
\end{minipage}\hfill
\begin{minipage}{0.49\textwidth}
\psfrag{t20}[c]{$t=20$}
\psfrag{u}[c]{$u(x,20)$}
\psfrag{x}[c]{$x$}
\includegraphics[width=\textwidth]{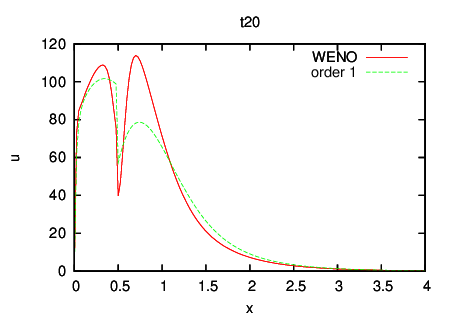}
\end{minipage}

\

\caption{Comparison between the WENO scheme and a first order scheme for the initial size distribution \eqref{eq:initial}, with the depolymerization value $k_\text{off}=8\times10^{-6}s^{-1}$ and the flux splitting parameter $\lambda=0.5$}\label{fig:1vs5}
\end{figure}

\begin{figure}[h]

\begin{minipage}{0.49\textwidth}
\psfrag{t0}[c]{$t=0$}
\psfrag{u}[c]{$u(0,x)$}
\psfrag{x}[c]{$x$}
\includegraphics[width=\textwidth]{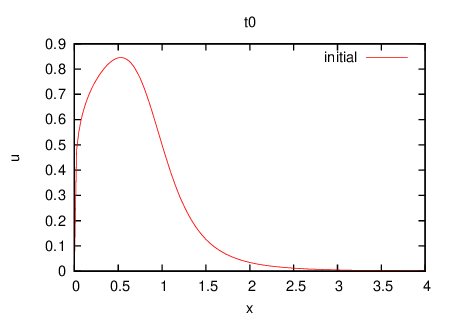}
\end{minipage}\hfill
\begin{minipage}{0.49\textwidth}
\psfrag{t1}[c]{$t=0.5$}
\psfrag{u}[c]{$u(x,0.5)$}
\psfrag{x}[c]{$x$}
\includegraphics[width=\textwidth]{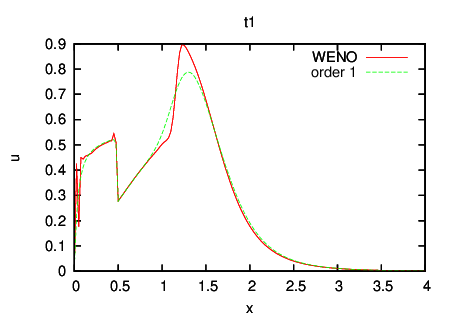}
\end{minipage}\\


\begin{minipage}{0.49\textwidth}
\psfrag{t12}[c]{$t=12$}
\psfrag{u}[c]{$u(x,12)$}
\psfrag{x}[c]{$x$}
\hspace{.1cm}\includegraphics[width=\textwidth]{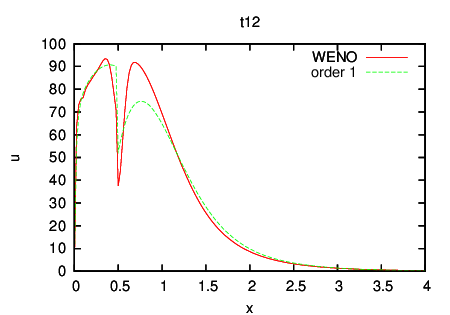}
\end{minipage}\hfill
\begin{minipage}{0.49\textwidth}
\psfrag{t20}[c]{$t=20$}
\psfrag{u}[c]{$u(x,20)$}
\psfrag{x}[c]{$x$}
\includegraphics[width=\textwidth]{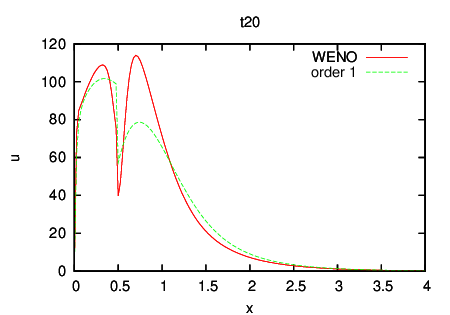}
\end{minipage}

\

\caption{Comparison between the WENO scheme and a first order scheme for a regular initial size distribution, with the depolymerization value $k_\text{off}=8\times10^{-6}s^{-1}$ and the flux splitting parameter $\lambda=0.5$}\label{fig:1vs5reg}
\end{figure}

We have considered that the mean size of the polymers at the initial time $t=0$ was $1380.$ This size can be multiplied by $5$ along the polymerization process (see Figure~\ref{fig:1vs5} keeping in mind that the mean size is represented by $0.2$ in this continuous model). So if we want to solve the discrete model, we have to consider a system of dimension close to $5000,$ and the computations are very heavy. If we limit this value to $200$ keeping the discrete model, then we lose a lot of precision. It is the same for the continuous model if it is discretized with a first order scheme. That is why we use a high order discretization, and we can see the difference in Figures~\ref{fig:1vs5}~and~\ref{fig:1vs5reg} : the high order scheme is able to capture strong variations of amplitude while the first order flattens them. We also remark that the size repartition converges to a bimodal distribution, as observed by \cite{Silveira}. This asymptotic profile is independent of the initial data (see the time $t=20$ in Figures~\ref{fig:1vs5}~and~\ref{fig:1vs5reg}) and can be seen as an eigenvector of the operator $\Q-\p_x\T$ (see \cite{M1,DG}).

\begin{figure}[h]

\begin{center}
\begin{minipage}{0.65\textwidth}
\psfrag{V}[l]{$V(t)$}
\psfrag{t}[c]{$t$}
\psfrag{n2}[c]{\small $\eta=2$}
\psfrag{n5}[c]{\small $\eta=5$}
\psfrag{n8}[c]{\small $\eta=8$}
\includegraphics[width=\textwidth]{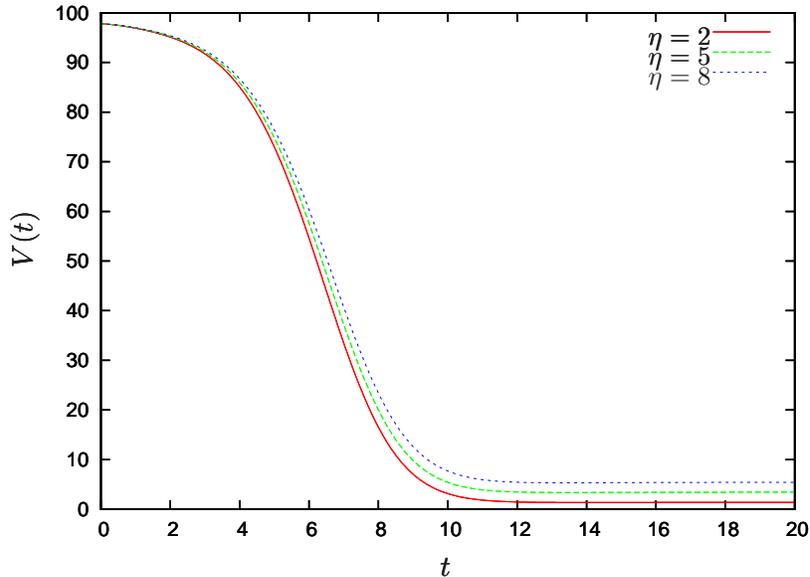}
\end{minipage}
\end{center}

\caption{Evolution of the quantity of monomers for different depolymerization rates, with the scheme $\G_{0.2}.$}\label{fig:monomers}

\end{figure}

The evolution of the quantity of monomers $V(t)$ is plotted in Figure~\ref{fig:monomers} for different values of the depolymerization rate $k_\text{off}.$ This quantity decreases since the monomers aggregate to polymers. Thus the mass of polymers increases and the speed of this evolution is similar to those observed by \cite{Knowles}. Concerning the dependence on $k_\text{off},$ the difference between the three curves is more significant when the time increases. For small times, when $V(t)$ is close to $V_0=98\mu M,$ the depolymerization can be neglected since $k_\text{off}$ is small compared to the product $V(t)k_\text{on}(x).$ Conversely, the equilibrium is reached when $\f{d}{dt}V(t)=0,$ so when $k_\text{off}\simeq Vk_\text{on}$ (see Equation~\eqref{eq:model}). That is why variations of $\eta$ influence essentially the ratio between the quantity of monomers and the mass of polymers at the equilibrium as we can see in Figure~\ref{fig:monomers}.

\section{Conclusion and future work}

\

We have written a high order conservative scheme for a polymerization-type equation. The choice of the flux splitting for the transport term has been discussed but the accurate decomposition remains unclear. It seems that unstabilities can be avoided by adapting the value of $\lambda$ but the oscillations remain present for any choice of the flux decomposition, even for a regular initial size distribution. A possible explanation for these phenomena can be that the integration method is not ``positive'' for the intervals at the boundaries. These points remain to be investigated for a better understanding and improvement.

As we have remarked in Section~\ref{sse:results}, the size distribution converges toward an equilibrium which corresponds to an eigenvector. The high-order WENO scheme presented in this paper could be used to numerically compute such eigenvectors. Another application of the code is to solve inverse problems (see \cite{DPZ,PZ}) in order to determine the size dependence of the different parameters.

\

\subsection*{Acknowledgement}

The authors are thankfull to Vincent Calvez, Marie Doumic Jauffret, Fr\'ed\'eric Lagouti\`ere and Natacha Lenuzza for their help and support during the Cemracs'09 research session.\\
This work has been done with the financial supports of\\
- the ANR contract TOPPAZ, allocation grant No. 4243, http://www-roc.inria.fr/bang/TOPPAZ/\\
- the CEA-Institute of Emerging Diseases and Innovative Therapies, Route du Panorama, Bat.60, F-92265 Fontenay-aux-Roses.

\bibliographystyle{abbrv}
\bibliography{Prion}

\end{document}